\documentclass[a4paper,12pt,legno]{article}
\usepackage[polish]{babel}
\usepackage[cp1250]{inputenc}
\usepackage[T1]{fontenc}
\usepackage{amsmath}
\usepackage{amsfonts}
\usepackage{amsthm}
\usepackage{mathrsfs}

\newtheorem{prop}[]{Proposition}
\newtheorem{lem}[]{Lemma}
\newtheorem{cor}{Corollary}
\newtheorem{ex}[]{Example}
\newtheorem{rem}[]{Remark}

\newcommand{\Rset}{\mathbb{R}}

\newcommand{\J}{\mathcal{J}}
\newcommand{\I}{\mathcal{I}}

\newcommand{\La}{\mathcal{L}}

\newcommand{\al}{\alpha}
\newcommand{\be}{\beta}

\newcommand{\bl}{\bigl(}
\newcommand{\br}{\bigr)}

\begin{document}

\title{Invariant measures under random integral mappings and marginal distributions of fractional L\'evy
processes.\footnote{Research financed by NCN no
DEC2011/01/B/ST1/01257}}

\author{Zbigniew J. Jurek}

\date{June 10, 2012.}

\maketitle
\begin{quote} \textbf{Abstract.} It is shown that some convolution semigroups of infinitely divisible
measures are invariant under the random integral mappings
$I^{h,r}_{(a,b]}$ defined in $(\star)$ below. The converse
implication is specified for the semigroups of generalized
s-selfdecomposable and selfdecomposable distributions. Some
application are given to the moving average fractional L\'evy
process (MAFLP).



\emph{Abbreviated title:} Invariant measures under random integral
mappings

\end{quote}

\medskip
Let us recall that \emph{the moving
average fractional L\'evy process} (in short, MAFLP) $(Z(t), t\in
\Rset)$ is given as follows
\begin{equation}
Z(t):= \int_{\Rset}\big((t-s)_+^{\al}-(-s)_+^{\al}\big)dY_{\nu}(s),
\ \ t\ge 0, \ \
\end{equation}
where $(Y_{\nu}(t), t \in \Rset)$ is a L\'evy process in $\Rset^d$,
$\nu$ is the probability distribution of the process at time $t=1$,
the parameter $\al$ is from the interval $(0,1/2)$ and $a_+ :
=\max(0, a)$ is the positive part of $a$. We study here how the laws
of $Z(t)$, in (1), is related  to  the law of $Y_{\nu}(1)$. Our
approach to that questions  is based on the so-called \emph{random
integral representation (or random integral mapping)}. This is a
technique that represents an infinitely divisible distribution, say
$\rho$, as a law of a random integral of the following form:
\begin{multline}
\qquad \rho= I^{h,r}_{(a,b]}(\nu):= \mathcal{L}\Big(\int_{(a,b]}\,
h(t)dY_{\nu}(r(t))\Big), \ \ \  \ \ \ \ (\star)   \\
\mbox{where} \ (a,b] \subset \Rset^+, \ h:\Rset^+ \to \Rset, \ \
Y_{\nu}(\cdot) \ \mbox{is a L\'evy process such that} \\
\mathcal{L}(Y_{\nu}(1))= \nu \ \mbox{and}\ r:\Rset^+ \to \Rset^+ \
\mbox{is a piecewise monotone time change},
\end{multline}
and its limit as $b\to \infty$;  cf. Jurek (2011) for a review of
the random integral mapping method and its application to
characterizations of classes of infinitely divisible laws. In this
context one might look at the conjectured \ "meta-theorem"  \ in \
\emph{The Conjecture}  on \ www.math.uni.wroc.pl/$\sim$zjjurek \, or
in Jurek (1985), p. 607  and Jurek (1988), p. 474.

We investigate classes of probability measures that are invariant
under random integral mappings $I^{h,r}_{(a,b]}$  (Proposition 1),
then we characterize those generalized s-selfedecomposable measures
that are, indeed, selfdecomposable ones (Proposition 2), and finally
we specify our results to the moving average fractional L\'evy
processes (MAFLP). [ Note the remark at the end of this paper.]

\textbf{1. Notations and the results.} Let $ID$ and $ID_{\log}$
denote the class of all infinitely divisible probability measures on
$\Rset^d$ and those that integrate the logarithmic function
$\log(1+||x||)$, respectively. Further, let  $\ast$ and
$\Rightarrow$ stand for  the convolution and the weak convergence of
measures, respectively. Thus $(ID,\, \ast, \, \Rightarrow)$ becomes
closed convolution subsemigroup of the semigroup of all probability
measures $\mathcal{P}$ (on $\Rset^d$).

Let $(Y_{\nu}(t), t\ge0)$ denotes a L\'evy process, i.e., a
stochastic process with stationary independent increments, starting
from zero, and with paths that continuous from the right and with
finite left limits (in short: cadlag), such that $\nu$ is its
probability distribution at time 1: $\mathcal{L}(Y_{\nu}(1))=\nu$,
where $\nu$ can be any $ID$ probability measure. Throughout the
paper $\mathcal{L}(X)$ will denote the probability distribution of
an $\Rset^d$-valued random vector $X$. Furthermore, for a
probability Borel measures $\mu$  its \emph{characteristic function}
$\hat{\mu}$ is defined as
\[
\hat{\mu}(y):=\int_{\Rset^d} e^{i<y,x>}\mu(dx), \ y\in\Rset^d,
\]
where $<\cdot,\cdot>$ denotes the scalar product (or a bilinear form
in case of  Banach space; cf. the concluding remark).

In  Section 3, formula (16), the L\'evy-Khintchine representation
for $\hat{\mu}$ of $\mu\in ID$ is recalled.

For the three parameters in (2) (i.e., the functions $h,r$ and the
interval $(a,b]$), let $\mathcal{D}^{h,r}_{(a,b]}$ denotes the
domain of definition of the mapping $I^{h,r}_{(a,b]}$. That is, the
set of all infinitely divisible measures $\nu$ ( L\'evy processes
$(Y_{\nu}(t), t\ge 0)$) such that the integral (2) is well defined.
Then the random integral mapping
\begin{equation}
I^{h,r}_{(a,b]}\ : \ \mathcal{D}^{h,r}_{(a,b]}\longrightarrow   ID,
\end{equation}
is a homomorphism between the corresponding convolution semigroups
because approximating $(\star)$ in (2) by the Riemann-Stieltjes sums
we get
\begin{equation}
\log\widehat{(I^{h,r}_{(a,b]}(\nu))}(y)=\int_{(a,b]}\log\widehat{\nu}(h(t)y)dr(t)
, \ y \in \Rset^d\,;
\end{equation}
 cf. for more details  Jurek-Vervaat (1983), Lemma 1.1 or Jurek and
Mason (1993), Chapter 3.
\begin{rem}
\emph{Cohen and Maejima (2011) defined the integral (1) in the same
way as it was in Marquardt (2006); see also the reference therein.
In particular they worked in the framework of L\'evy processes with
finite variance, square integrable functions and Euclidean spaces.
However, using the formal integration by parts we are able to define
random integrals for larger class of integrands $h$ and L\'evy
processes $Y$. Moreover, still having the crucial equality (4).}
\end{rem}

From our definition of random integrals, in particular from (4), we
infer the following properties:
\begin{equation}
I^{h,r}_{(a,b]}(\nu_1)\ast I^{h,r}_{(a,b]}(\nu_2) =
I^{h,r}_{(a,b]}(\nu_1 \ast \nu_2), \ \ \
I^{h,r}_{(a,b]}(T_u\nu)=I^{u h,r}_{(a,b]}(\nu)=T_u\big(
I^{h,r}_{(a,b]}(\nu)\big)
\end{equation}
\begin{equation}
I^{h,r}_{(a,b]}(\nu^{\ast s})= (I^{h,r}_{(a,b]}(\nu))^{\ast s} =
(I^{h,\,sr}_{(a,b]}(\nu)),\
I^{h,r}_{(a,b]\cup(b,c]}(\nu)=I^{h,r}_{(a,b]}(\nu)\ast
I^{h,r}_{(b,c]}(\nu)
\end{equation}
\begin{equation}
\mbox{if} \ \ \nu_n\Rightarrow\nu \ \ \mbox{then} \ \ I^{h,r}_{(a,b]}(\nu_n)
\Rightarrow I^{h,r}_{(a,b]}(\nu),
\end{equation}
where $T_u$ is the dilation, i.e., $T_u(x):=ux$, $u \in \Rset$,
$x\in\Rset^d$ and $s>0$. [ Replacing the dilation $T_u$ by a matrix
( or a bounded linear operator) $A$ in (5) we get
$A(I^{h,r}_{(a,b]}(\nu))=I^{h,r}_{(a,b]}(A\nu)$.]

Random integrals over half-lines are defined as limits almost surely
(or in distribution or in probability) over finite intervals $(a,b]$
as $b\to \infty$; cf. Jurek and Vervaat (1983).

\begin{prop}
Let $\mathcal{K}$ be a closed convolution subsemigroup of the
semigroup $ID$ (of all infinitely divisible measures ) that is also
closed under dilations  and convolution powers (i.e., if $a\in
\Rset$ and $\nu\in\mathcal{K}$ then $T_a\nu \in \mathcal{K}$ and for
$c>0$  also $\nu^{\ast c}\in \mathcal{K}$). Then if $\nu\in
\mathcal{K}\cap \mathcal{D}^{h,r}_{(a,b]}$ then
$I^{h,r}_{(a,b]}(\nu) \in \mathcal{K}.$

The same holds also for improper random integrals (over half-lines or lines) provided they are well-defined.
\end{prop}
Using the properties of (5)-(7) we get
\begin{cor}
Domains of definition $\mathcal{D}^{h,r}_{(a,b]}$ of random
integrals $I^{h,r}_{(a,b]}$ are examples of semigroups $\mathcal{K}$
from Proposition 1.
\end{cor}
Other, more explicite, examples of classes $\mathcal{K}$ are given
in Example 1 below, after introducing some auxiliary notions and
notations.

For the purpose of this note we will consider two specific random
integral mappings and their  corresponding semigroups.

Firstly, for $\beta>0$ and $\nu\in
ID$, let us define
\begin{equation}
I^{t,\,t^{\beta}}_{(0,1]}(\nu)\equiv\mathcal{J}^{\beta}(\nu):\,=\La\bl\int_{(0,1]}
t\;dY_{\nu}(t^{\beta})\br, \ \ \mbox{and} \ \
\mathcal{U}_{\beta}:\,=\mathcal{J}^{\beta}(ID).
\end{equation}
To the distributions from the semigroups $\mathcal{U}_{\beta}$ we
refer to as\textit{ the generalized s-selfdecomposable distributions}.
\begin{rem}
\emph{The classes $\mathcal{U}_{\beta}$ were already introduced in
Jurek (1988) as the limiting distributions in some schemes of
summing of independent variables. The terminology has its origin in
the fact that distributions from the class $\mathcal{U}_{1} \equiv
\mathcal{U}$ were called \emph{s-selfdecomposable distribution} (the
"\emph{s}-" , stands here for \emph{the shrinking operations} that
were used originally in the definition of $\mathcal{U}$); cf. Jurek
(1981), (1985), (1988) and references therein.}
\end{rem}
Secondly, for $\nu\in ID_{\log}$ let us put
\begin{equation}
I^{e^{-t},\,t}_{(0,\infty)}(\nu)\equiv \mathcal{I}(\nu):=
\mathcal{L}\bigl(\int_{(0,\infty)}e^{-s}\, d\,Y_{\nu}(s)\bigr) \ \
\mbox{and} \ \ L:=\mathcal{I}(ID_{\log}).
\end{equation}
The distributions from the semigroup $L$ are called
\emph{selfdecomposable} ones or \emph{L\'evy class L distributions}.
Let us stress here that the logarithmic moment guarantees the
existence of the improper random integral (7); cf. Jurek-Vervaat
(1983), Theorem 2.3 or Jurek-Mason (1993, Chapter III.
\begin{rem}
\emph{In classical probability  theory the selfdecomposability  is
usually defined via some decomposability property or by scheme of
limiting distributions. However,  since Jurek-Vervaat (1983) we know
that the class $L$ coincides with the class of distributions of
random integrals given in (9). Hence it is used here as its
definition.}
\end{rem}
Between the classes $L$, $\mathcal{U}_{\beta} \ (\be>0)$, the class
$\mathcal{G}$ of all Gaussian measures and the class $\mathcal{S}$
of all stable probability measures we have the following proper
inclusions:
\begin{equation}
\mathcal{G}\subset \mathcal{S} \subset L \subset \mathcal{U}_{\beta}
\subset ID, \ \ \mbox{i.e.,} \ \ \mathcal{I}(ID_{\log}) \subset
\mathcal{J}^{\beta}(ID)\,.
\end{equation}
\begin{rem}
\emph{It might be of an interest to recall here that many classical
distributions in mathematical statistics such as gamma, t-Student,
Fisher F etc. are in the class L but, of course, they are not
stable; cf. the survey article Jurek (1997) or Jurek-Yor (2004) or
the book by Bondesson (1992).}
\end{rem}

\begin{ex}
The classes $L$ (of the selfdecomposable distributions),
$\mathcal{U}_{\beta}$ (of the generalized s-selfdecomposable
distributions) and $\mathcal{G}$ (of the Gaussian measures) are
examples of the above class $\mathcal{K}$. Also the Urbanik class
$L_{\infty}$, that coincides with the smallest closed convolution
semigroup generated by all stable distributions, is an example of
the class  $\mathcal{K}$; cf. Urbanik (1973), or Jurek (2004).
\end{ex}
From the inclusions in (10) we get that all selfdecomposable
measures are generalized s-selfdecomposable ones whenever $\be>0$.
With the notations described below  the formula (16), we give
conditions for the converse claim.
\begin{prop}
Let  $\nu =[b, S, N]\in ID$ and $\rho =[a,R,M]\in ID_{\log}$. Then
the following conditions are equivalent:
\item[(i)] $\J^\be(\nu)= \I(\rho)$, i.e., a generalized
s-selfdecomposable measure is in fact a selfdecomposable one;
\item[(ii)] $\nu=\rho^{\ast 1/\beta} \ast\mathcal{I}(\rho)$ and,
consequently, $\nu\in ID_{\log}$;
\item [(iii)] $\Rset^d \ni y \to
\exp\be\big[\log\hat{\nu}(y)-\be\int_0^1\log\hat{\nu}(t\,y) t^{\be-1}dt
\big] \ \mbox{is a Fourier transform}$ \\  \ \ \ \ \ \ \ \ \
\mbox{of an $ID_{\log}$ measure};
\item[(iv)] $\int_0^1\big(N(A)-N(s^{-1}A)\big)s^{\be-1}ds\ge 0 \ \
\mbox{for all Borel sets}\ A \ \mbox{such that} \ 0 \notin A $ and
$\int_{(||x||>1)}\log||x||N(dx)<\infty$.
\end{prop}
Here we have the above condition (i) in terms of the triples from
L\'evy -Khintchine representation :
\begin{cor}
In order that $\J^\be([b,S,N])=\I([a,R,M])$ it is necessary and
sufficient that
\begin{multline*}
b=(\be+1)\be^{-1}\,a+\int_{(||x||>1)} x\,||x||^{-1}M(dx) \ \
\mbox{and} \  \ S=(\be+2)(2\be)^{-1}\,R \\
\ \  \  \mbox{and} \ \ \ N(A)=\be^{-1}M(A)
+\int_0^1M(t^{-1}A)t^{-1}dt\ \ \mbox{for all}  \ \ A \in
\mathcal{B}_0 .
\end{multline*}
\end{cor}

\textbf{2. The case of MAFLP.} Now we will specify our
considerations to the case of MAFLP $Z(t)$ given in (1). First of
all, note that similarly as in (1), for a L\'evy process
$(Y_{\nu}(t), t\ge 0)$, putting
\begin{equation}
U^{(\nu)}(t):=\int_{-\infty}^0((t-s)^{\al}-(-s)^{\al})dY_{\nu}(s) \
\mbox{and}\ V^{(\nu)}(t):=\int_0^t(t-s)^{\al}dY_{\nu}(s)
\end{equation}
we get that
\begin{equation}
Z(t)=U^{(\nu)}(t)+V^{(\nu)}(t) \ \ \mbox{and the summands are
independent.}
\end{equation}
This is so because two-sided L\'evy process (i.e., with time index
in $\Rset$) is defined by taking independent copies of L\'evy
processes on both half-lines; cf. Marquardt (2006), p. 1102.

\noindent Furthermore, using the invariance principle for L\'evy
processes, that is, the property that for each fixed positive t we
have
\[
\big(-Y_{\nu}(-s), 0\le s \le t\big)\stackrel{d}{=} \big(Y_{\nu}(s),
0 \le s \le t\big)\stackrel{d}{=} \big(Y_{\nu}(t)-Y_{\nu}(t-s)-,
0\le s\le t\big)
\]
(the equality in distribution of three L\'evy processes) we infer
that that
\begin{equation}
U^{(\nu)}(t)\stackrel{d}{=}\int_0^{\infty}\big((t+s)^{\al}-s^{\al}\big)dY_{\nu}(s),
\ \ \ \  \ \
V^{(\nu)}(t)\stackrel{d}{=}\int_0^t\,s^{\alpha}dY_{\nu}(s).
\end{equation}
Of course, from (13) and (8)  we have that
\[
I^{s,\,s^{1/\alpha}}_{(0,1]}(\nu)=I^{s^{\al},s}_{(0,1]}(\nu)=
\mathcal{L}(V^{(\nu)}(1))\in \mathcal{U}_{1/\al} \ \ \ \mbox{and} \
\ 2<1/\al.
\]
Then for $t>0$, the above with (5), (6) and Example 1 (for the class
$\mathcal{U}_{\be}$) give
\begin{multline}
T_{t^{\al}}\big[\big(I^{\,s^{\al}, \,s}_{(0,1]}(\nu) \big)^{\ast
t}\big]=T_{t^{\al}}\big[I^{\,s^{\al}, \,t\,s}_{(0,1]}(\nu)\big]\\
=I^{\,(t\,s)^{\al}, \,t\,s}_{(0,1]}(\nu) =I^{\,s^{\al},
\,s}_{(0,t]}(\nu)=\mathcal{L}(V^{(\nu)}(t))\in \mathcal{U}_{1/\al},
\qquad
\end{multline}
and consequently we get
\begin{cor}
For all infinitely divisible measures $\nu$, probability
distributions of $V^{(\nu)}(t)$ are in the class
$\mathcal{U}_{1/\al}$ of generalized s-selfdecomposable probability
measures with $1/\al>2$.
\end{cor}

In (13) integrals $U^{(\nu)}(t)$ over half-line are defined as
limits, i.e.,
\begin{equation}
U^{(\nu)}(t)=  \lim_{b \to \infty}U^{(\nu), b}(t):=\lim_{b\to
\infty} \int_{(0,b]}\big((t+s)^{\al}-s^{\al}\big)dY_{\nu}(s)
\end{equation}
a.s. or in distribution. Because of (11) and (13),
$U^{(\nu),\,b}(t)$ and $V^{(\nu)}(t)$ are stochastically independent
and $\lim_{b\to\infty}[U^{(\nu), b}(t)+V^{(\nu)}(t)]=Z(t)$.

Since  an integral $U^{(\nu),\,b}(t)$ is of the form
$I^{h,r}_{(a,b]}$ we may apply Proposition 1 and get properties of
marginal distributions of MAFLP summarized as follows:
\begin{cor}
Let $\mathcal{K}$ be a closed convolution semigroup of infinitely
divisible measures that is also closed under dilations  and
convolution powers (i.e., if $c>0$ and $\nu\in\mathcal{K}$ then
$T_c\nu \in \mathcal{K}$ and $\nu^{\ast c}\in \mathcal{K}$). Then

(a) if $\nu\in \mathcal{K}$ then
$\mathcal{L}[U^{(\nu),b}(t)+V^{(\nu)}(t)] \in \mathcal{K}$ for all
$t>0$;

(b) if $\nu\in \mathcal{K}$ and  MAFLP $Z(.)$ is well defined then
its marginal distributions $\mathcal{L}(Z(t))\in \mathcal{K}$ for
all $t>0$.
\end{cor}

\begin{rem}
\emph{ The above corollary (part (b)) for  the case of
selfdecomposable measures was also noted in Cohen and Maejima
(2011).}
\end{rem}

\medskip
\textbf{3. Proofs.}  Recall that for infinitely divisible measures
$\mu$ their characteristic functions admit the following
L\'evy-Khintchine representation:
\begin{multline}
\hat{\mu}(y)= e^{\Phi(y)}, \ y \in \Rset^d, \ \ \mbox{and the L\'evy
exponents $\Phi$ are of the form} \\  \Phi(y)=i<y,a>-
\frac{1}{2}<y,Ry>  +  \qquad \qquad \qquad \\ \int_{\Rset^d
\backslash \{0 \}}[e^{i<y,x>}-1-i<y,x>1_B(x)]M(dx),
\end{multline}
where $a$ is  a \emph{shift vector}, $R$ is a \emph{covariance
operator} corresponding to the Gaussian part of $\mu$, $1_B$ is the
indicator function of the unit ball $B$ and $M$ is a \emph{L\'evy
spectral measure}. Since there is a one-to-one correspondence
between a measure $\mu \in ID$ and the triple $a$, $R$ and $M$ in
its L\'evy-Khintchine formula (10) we will write $\mu=[a,R,M]$; also
cf. the remark at the end of this paper.

\medskip
From (13) and an appropriate version of (4) we get the following
\begin{lem}
(i) If $\nu = [b,S,N]$ and
$\J^\be(\nu)=[b^{(\be)},S^{(\be)},N^{(\be)}]$ then
\begin{multline*}
b^{(\be)}:=
\frac{\be}{\be+1}\,(b+\int_{(||x||>1)}x\,||x||^{-1-\be}N(dx)\,); \ \
\
S^{(\be)}:=\tfrac{\be}{2+\be}\,S;  \qquad    \\
N^{(\be)}(A):=\int_0^1 \, N(t^{-1/\be}A) \,dt\,= \be\,\int_0^1
N(s^{-1}A)s^{\be-1}ds \ \ \mbox{for each} \,\,A \in \mathcal{B}_0
\end{multline*}
(ii) If $\mu = [a,R,M]$ and
$\mathcal{I}(\nu)=[a^{\sim},R^{\sim},M^{\sim}]$ then
\begin{multline*}
a^{\sim}:= a+\int_{(||x||>1)}x\,||x||^{-1}M(dx); \ \ \ R^{\sim}:= \frac{1}{2}\,R \\
M^{\sim}(A):=\int_0^{\infty}M(e^t\,A)dt=\int_0^1M(t^{-1}A)\,t^{-1}dt \ \ \ \mbox{for each} \,\,A \in\mathcal{B}_0; \ \ \ \
\end{multline*}
\end{lem}
\noindent cf. Czyżewska-Jankowska and Jurek (2011), Lemma 2, and
Jurek and Vervat (1983) on p. 250 for more details.

\medskip
\textit{Proof of Proposition 1 .} For $h$ of bounded variation,
cadlag L\'evy process $Y$ and montone $r$ we define here the random
integral $\star$ as follows:
\begin{multline}
\int_{(a,b]} h(t)dY(r(t)):= h(b)Y(r(b))-h(a)Y(r(a))-\int_{(a,b]}Y(r(t)-)dh(t)\\
=h(b)(Y(b)-Y(a))-\int_{(a,b]}(Y(r(t)-)-Y(r(a))dh(t),
\end{multline}
where  $Y(r(t)-)$ denotes the left-hand limit. Consequently, for the
partition $a=t_0<t_1<t_2<...<t_n=b$, the random integral (17) can be
approximated by the Riemann-Stieltjes sums
\begin{multline}
h(b)[Y(r(b))-h(a)Y(r(a))]-\sum_{j=1}^nY(r(t_j)(h(t_j)-h(t_{j-1}))\\
=\sum_{j=1}^nh(t_j)(Y(r(t_j)-Y(r(t_{j-1})). \ \ \ \  \ \ \ \ \ \
\end{multline}
The summands in (18) are independent and since
$\nu=\mathcal{L}(Y(1))\in\mathcal{K}$ we get that
\[
\mathcal{L}[h(t_j)(Y(r(t_j)-Y(r(t_{j-1}))]=T_{h(t_j)}(\mathcal{L}(Y(1))^{\ast
(r(t_j)-r(t_{j-1}))})\in \mathcal{K},
\]
if $r(t_j)-r(t_{j-1})\ge 0$. Similarly,
$\mathcal{L}[-h(t_j)(Y(r(t_{j-1})-Y(r(t_{j}))] \in\mathcal{K}$ when
$r(t_j)-r(t_{j-1})\le 0$. Closenesses and semigroup property of
$\mathcal{K}$ guarantees that $I^{h,r}_{(a,b]}(\nu)\in \mathcal{K}$,
which gives the proof of Proposition 1 for finite intervals $(a,b]$.
Since integrals on half-lines are given as weak limits of those over
$(a,b]$ as $b\to\infty$ and $\mathcal{K}$ is closed in weak topology
we get Proposition 1 for half-lines, which completes a proof.

\medskip
\textit{Proof of Proposition 2 .} $(i)\equiv (ii)$. Let us put
$\Phi(y):=\hat{\nu}(y)$ and $\Psi(y):=\hat{\rho}(y)$, i.e., they are
the corresponding L\'evy exponents. Then using (4), (8) and (9) we
infer that (i) is equivalent the following identity
\begin{equation}
\beta\int_0^1\Phi(ty)t^{\beta-1}dt=\int_0^{\infty}\Psi(e^{-s}y)ds=\int_0^1\Psi(ty)\frac{dt}{t},
\ \ \ \ \mbox{for all y}\in \Rset^d.
\end{equation}
Putting into above $sy$  for $y$, where $s\in\Rset$ varies and y is fixed, and then substituting $w:=st$ we get
\begin{equation*}
\int_0^s\Phi(wy)w^{\beta-1}dw=\be^{-1}s^{\be}\int_0^s\Psi(wy)\frac{dw}{w}.
\end{equation*}
Differentiating with respect to $s$ and then putting $s=1$ we arrive  at
\begin{equation}
\Phi(y)=\int_0^1\Psi(wy)\frac{dw}{w} +\be^{-1}\Psi(y), \ \ \
\mbox{for all y},
\end{equation}
and after exponentiating both sides we get the equality (ii) in
terms of Fourier transforms.

Conversely, starting with (20) and substituting  $ty$ for $y$ and
then integrating both sides over the unit interval with respect
$dt^{\be}$ we arrive at
\[
\beta\int_0^1\Phi(ty)t^{\beta-1}dt=\int_0^1\Psi(ty)\frac{dt}{t}
\]
which means that $\J^\be(\nu)= \I(\rho)$.

$[(i)\equiv (ii)]\Rightarrow (iii)$. Substituting $\mathcal{J}^{\beta}(\nu)$ for $\mathcal{I}(\rho)$ in (ii)
and then taking Fourier transforms both sides we get that $\hat{\rho}$ has the form as in (iii).

$(iii)\Rightarrow (iv)$. Let $\rho\in ID_{\log}$ has Fourier
transform given by (iii). Then $\rho^{\ast 1/\be}\in ID_{\log}$ and
its L\'evy spectral measure is of the form
\[
N(A)-N^{(\be)}(A)=\be\int_0^1(N(A)-N(t^{-1}A)\,t^{\be-1}dt\ge 0 \ \
\mbox{for all} \ \ A\in\mathcal{B}_0
\]
which is  the claim (iv).

$(iv)\Rightarrow (i).$  Multiplying (iv) by $\be$, and using the
notation from Lemma 1 (i), we have that $0\le N^{(\be)}\le N$.
Consequently, $N-N^{(\be)}$ is a L\'evy spectral measure with finite
log-moment; cf. Czyżewska-Jankowska and Jurek (2011), Lemma 2 (ii).
Furthermore, from Lemma 1(ii),
\begin{multline*}
(\be (N-N^{(\be)}))^{\sim}(A)
 = \be(\int_0^1N(t^{-1}A)t^{-1}\,dt-
\int_0^1\int_0^1\,N(t^{-1/\be}s^{-1}A)dt\,s^{-1}ds) \\
=
\be(\int_0^1N(t^{-1}A)t^{-1}\,dt-\int_0^1(\int_0^sN(w^{-1}A)\be w^{\be-1}dw)s^{-(\be+1)}ds)\\
=\be(\int_0^1N(t^{-1}A)t^{-1}\,dt-\int_0^1N(w^{-1}A)
w^{\be-1}(w^{-\be}-1)dw)\\
= \be\int_0^1N(w^{-1}A)w^{\be-1}dw= N^{(\be)}(A). \ \qquad \qquad
\end{multline*}
Similarly, using Lemma 1 (ii), for Gaussian covariance operator we
have
\[
(\be(S-S^{(\be)}))^{\sim}=
(2\be(\be+2)^{-1}S)^{\sim}=\be(\be+2)^{-1}S=S^{(\be)}.
\]
Finally,  applying (ii) in Lemma 1 for the shift vectors we get
\begin{multline*}
(\be(b-b^{(\be)}))^{\sim}=\be(b-b^{(\be)})+\be\int_{(||x||>1)}\frac{x}{||x||}N(dx)-\be\int_{(||x||>1)}\frac{x}{||x||}N^{(\be)}(dx)\\
=\be(b-b^{(\be)})+\be\int_{(||x||>1)}\frac{x}{||x||}N(dx)-\be^2\int_0^1\int_{(||x||>1)}\frac{x}{||x||}N(s^{-1}dx)s^{\be-1}ds \\
=\be(b-b^{(\be)})+\be\int_{(||x||>1)}\frac{x}{||x||}N(dx)-
\be\int_{(||w||>1)}\frac{w}{||w||}\int_{||w||^{-1}}^1\be
s^{\be-1}N(dw)\\
=\be(b-b^{(\be)})+\be\int_{(||x||>1)}\frac{x}{||x||}N(dx)-\be[\int_{(||w||>1)}\frac{w}{||w||}N(dw)-\int_{(||w||>1)}\frac{w}{||w||^{1+\be}}N(dw)\\
=\be(b+\int_{(||w||>1)}\frac{w}{||w||^{1+\be}}N(dw)\,-b^{(\be)})=\be
((\be+1)\be^{-1}b^{(\be)}-b^{(\be)})=b^{(\be)}.
\end{multline*}
All in all we have that
$\rho=[\,\be(b-b^{(\be)}),\be(S-S^{(\be)}),\be (N-N^{(\be)})\,]\in
ID_{\log}$ and $\mathcal{I}(\rho)=\mathcal{J}^{\beta}(\nu)$, which
completes the proof of $(iv)\Rightarrow (i)$ and thus the proof of
Proposition 2.

\medskip
\textbf{Concluding remark.} Last but not least, although we
presented our results for the Euclidean space $\Rset^d$,  our
methods are applicable for processes and random variables with
values in any real separable Banach space. For an exposition of
probability on Banach spaces see Araujo-Gin\'e (1980) or
Ledoux-Talagrand (1991) and for a case of Hilbert spaces we
recommend  Parthasarathy (1967), Chapter VI. In particular, the
crucial L\'evy-Khintchine representation (16) holds true in the
generality of separable infinite dimensional Banach spaces. But
there is no integrability criterium for L\'evy (spectral) measures
$M$ analogous that we have on Euclidean and Hilbert spaces.

\medskip
\begin{center}
\textbf{References}
\end{center}

\medskip
\noindent [1] A. Araujo and E. Gin\'e (1980), \emph{The central
limit theorem for real and Banach valued random variables}, John
Wiley \& Sons, New York.

\noindent [2] L. Bondesson (1992), \emph{Generalized gamma
conolutions and related classes of distributions and densities},
Lecture Notes in Stat. vol. 76, Springer Verlag, New York.

\noindent [3] S. Cohen and M. Maejima (2011), Selfdecomposability of
moving average fractional L\'evy processes,  \emph{Stat.\& Probab.
Lett.} vol.81, pp. 1664-1669.

\noindent[4] A. Czyżewska-Jankowska and Z. J. Jurek (2011),
Factorization property of generalized s-selfdecompsable measures and
class $L^f$ distributions, \emph{Theory Probab. Appl.} vol. 55, No
4, pp. 692-698.

\noindent [5] Z. J. Jurek (1981), Limit distributions for sums of
shrunken random variables, \emph{Dissertationes Mathematicae}, vol.
\textbf{185}, 46 pp.

\noindent [6] Z. J. Jurek (1985). Relations between the
s-selfdecomposable and selfdecomposable measures. \emph{Ann.
Probab.} vol.13, Nr 2, str. 592-608.

\noindent [7] Z. J. Jurek (1988). Random Integral representation for
Classes of Limit Distributions Similar to L\'evy Class $L_{0}$,
 \emph{Probab. Th. Fields.} 78, pp. 473-490.

\noindent[8] Z. J. Jurek (1997). Selfdecomposability: an exception
or a rule ? , \emph{Ann. Univ. Marie Curie Skłodowska; Sec. A,
Math.}, \textbf{LI}.1, pp. 93-107

\noindent [9] Z. J. Jurek (2004). The random integral representation
hypothesis revisited: new classes of s-selfdecomposable laws.  In:
ABSTRACT AND APPLIED ANALYSIS; \emph{Proc. International Conf.
Hanoi, Vietnam 13-17 August 2002;} pp. 495-514. WORLD SCIENTIFIC,
June 2004.

\noindent[10] Z. J. Jurek (2011). The random integral representation
conjecture: a quarter of a century later. \emph{Lithuanian Math.
Journal} vol. no 3 pp.

\noindent[11] Z. J. Jurek and J.D. Mason (1993).
\emph{Operator-limit distributions in probability}, J. Wiley $\&$
Sons,  New York.

 \noindent [12] Z. J. Jurek and W. Vervaat (1983). An
integral representation for selfdecomposable Banach space valued
random variables, \emph{Z. Wahrscheinlichkeitstheorie verw.
Gebiete}, 62, pp. 247-262.

\noindent [13] Z. J. Jurek and M. Yor (2004). Selfdecomposable laws
associated with hyperbolic functions, \emph{Probab. Math. Stat.} 24,
no.1, pp. 180-190.

\noindent (www.math.uni.wroc.pl/$\sim$pms)

\noindent[14] M. Ledoux and M. Talagrand (1991).  \emph{Probability
in Banach spaces}, Springer-Verlag.

\noindent[15] T. Marquardt (2006). Fractional L\'evy processes with
application to long memory moving average
processes,\emph{Bernoulli}, 12 no 6, pp.1099-1126.

\noindent[16] K. R. Parthasarathy (1967). \emph{Probability measures
on metric spaces}. Academic Press, New York and London.

\noindent [17] K. Urbanik (1973). Limit laws for sequences of normed
sums satisfying some stability conditions. In: \emph{Multivariate
Analysis III} (P.R. Krishnaiah, Ed.), pp. 225 -237, Academic Press,
New York.

\medskip
\noindent
Institute of Mathematics \\
University of Wroc\l aw \\
Pl.Grunwaldzki 2/4 \\
50-384 Wroc\l aw, Poland \\
e-mail: zjjurek@math.uni.wroc.pl \ \ \ \  \
www.math.uni.wroc.pl/$^{\sim}$zjjurek

\end{document}